\documentclass[B,12pt]{amsart}
\usepackage{graphicx}
\newtheorem{thm}{Theorem}[section]

\newtheorem{lem}[thm]{Lemma}

\theoremstyle{definition}
\newtheorem{defn}[thm]{Definition}
\theoremstyle{remark}
\newtheorem{rem}[thm]{Remark}
\numberwithin{equation}{section}

\begin{document}
\title[Bounded and Almost Periodic Solutions]{Bounded and Almost Periodic Solutions for Second Order Differential Equation
Involving Reflection of the Argument}%
\author{Daxiong Piao $^\dag$ $^\sharp$ }
\author{Na Xin $^\ddag$}
\address{$^{\dag}$School of Mathematical Sciences, Ocean University of China,
Qingdao, 266071,China}
\address{${^\ddag}$Department of Mathematics,Normal College, Qingdao
University, Qingdao, 266071, China}
\email{dxpiao@ouc.edu.cn (D.Piao)  \,\, xinna1026@sina.com(Na Xin)}%
\thanks{$^\sharp$ Corresponding author. Supported by the NSF of Shangdong
Province (No.ZR2012AM018)}


\begin{abstract}
In this paper we investigate the existence and uniqueness of
bounded, periodic and almost periodic solutions for second order
differential equations involving reflection of the argument.The
relationship between frequency modules of forced term and the
solution of the equation is considered.
\end{abstract}
\maketitle

\section{Introduction}

The differential equations involving reflection of argument have
applications in the study of stability of differential-difference
equations, see Sarkovskii [1] , and such equations show very
interesting properties by themselves, so many authors worked on
them. Wiener and Aftabizadeh [2] initiated to study boundary value
problems for the second order differential equations involving
reflection of the argument. Gupta [3,4] investigated two point
boundary value problems for this kind of equations under the
Caratheodory conditions. Afidabizadeh, Huang, and Wiener [5]
studied the existence of unique bounded solution of first order
equation
$$
  \dot{x}(t)=f(t,x(t),x(-t)).
$$
They proved that $x(t)$ is almost periodic by assuming the
existence of bounded solution $x(t)$ of the equation. In [6,7],
One of the present authors investigated existence and uniqueness
of periodic, almost periodic and pseudo almost periodic solutions
of the equations
$$
\dot {x}(t)+ax(t)+bx(-t)=g(t),b\neq 0, t\in R,
$$
and
$$
 \dot {x}(t)+ax(t)+bx(-t)=f(t,x(t),x(-t)), b\neq 0, t\in R.
$$

In [12], the authors studied the first order operator
$x'(t)+mx(-t)$ coupled with periodic boundary value conditions,and
described the eigenvalues of the operator and obtained the
expression of its related Green¡¯s function in the non resonant
case.

 Our present paper is motivated by above references, and devoted to
investigate the existence and uniqueness of bounded, periodic and
almost periodic solution of the second order linear equation
\begin{equation}\label{1.1}
 \ddot {x}(t)+ax(t)+bx(-t)=g(t), b\neq 0, t\in R,
\end{equation}
and nonlinear equation
\begin{equation}\label{1.2}
 \ddot {x}(t)+ax(t)+bx(-t)=f(t,x(t),x(-t)), b\neq 0, t\in R,
\end{equation}
respectively, where $g$ and $f$ satisfy some assumptions which
will be stated later.

Our paper is organized as following.In Section 2,we state some
lemmas and basic formulations; In Section 3,we study the case $a <
b < -a$ for the equations (1.1) and (1.2); In Section 4, we study
the case $-a < b < a $.

\medskip
\section{ Some Lemmas and useful formulations }

Now we give some definitions for our business.

{\defn $^{[8-9]}$ \,\,  A function $f:{\bf R}\rightarrow{\bf C}$ is
almost periodic , if the $\epsilon$--translation of $f$
$$
{\bf T}(f,\epsilon)=\{ \tau\in {\bf R}: |f(t+\tau)-f(t)|<\epsilon
,\forall t \in {\bf R} \}
$$
is relatively dense in $ {\bf R}$. We denote the set all such
functions by $AP({\bf R})$.

{\defn $^{[8-9]}$ \,\, A function $F:{\bf R}\times {\bf C
}^2\rightarrow R$ is almost periodic for $t$ uniformly on ${\bf C
}^2$ , if for any compact $ {\bf W} \subset  {\bf C}^2$, the
$\epsilon$--translation of $F$
$$ {\bf T}(F,\epsilon,{\bf W})=\{\tau\in{\bf
R}:|F(t+\tau,x,y)-F(t,x,y)|<\epsilon, \forall(t,x,y)\in {\bf
R}\times {\bf W} \}
$$
is relatively dense in ${\bf R}$ .} We denote the set of all such
functions by $AP({\bf R}\times {\bf C}^2)$. ${\bf R}$ is the set of
all real numbers ,and ${\bf C}$ is the set of all complex numbers.

we state two useful lemmas those can be easily proven.

{\lem  If $g(t)\in AP({\bf R})$, then $g(-t)\in AP({\bf R})$.
Furthermore if $\tau$ is an $\epsilon$-translation of $g(t)$, then
$\tau$ is also an $\epsilon$-translation of  $g(-t)$. If $g(t)$ is
$\omega$-periodic, then $g(-t)$ is also $\omega$-periodic.}

{\lem If $g(t)\in AP({\bf R})$, then $mod (g(t))= mod (g(-t))$ and
$Freq(g(t))=-Freq(g(-t))$, where $Freq(g)$ denotes the frequency
set of $g(t)$.}

We refer the readers to good books [8-10] for the  basic results
on the almost periodic functions .

Before treating  the nonlinear equation (1.2), we need consider
the linear equation (1.1) first .

    Let $x_1=x(t)$, $x_2=x_1(-t)$, $x_3=\dot x_1(t)$, $x_4=x_3(-t)$, then we obtain a system
\begin{equation}\label{2.1}
\left\{
\begin{array}{l}
\dot x_1=x_3\\
\dot x_2=-x_4\\
\dot x_3=-ax_1-bx_2+g(t)\\
\dot x_4= bx_1+ax_2-g(-t)
\end{array}\right.
\end{equation}

or

\begin{equation}\label{2.2}
\dot{\bf x}(t)=A {\bf x}+{\bf g}(t)
\end{equation} where
${\bf x}=\left(
\begin{array}{c}
x_1 \\
x_2\\
x_3\\
x_4\end{array}\right)$, $ A=\left(\begin{array}{cccc}
0&0 &1& 0\\
0&0 &0&-1\\
-a&-b&0&0\\
b&a& 0 &0
\end{array} \right) $ and ${\bf
g}(t)=\left(\begin{array}{c}
0\\
0\\
g(t)\\
-g(-t)
\end{array}
\right)$.

\medskip

 Now $det{(r I-A)}=r^4+2a r^2+a^2-b^2$, so the
eigenvalues of $A$ are $r_1=\alpha=\sqrt{b-a}$,
$r_2=-\alpha=-\sqrt{b-a},r_3=\beta=\sqrt{-a-b}$,
$r_4=-\beta=-\sqrt{-a-b}$. If $\alpha\beta\neq 0$,then their
corresponding eigenvectors are

${\bf v}_1=\left(\begin{array}{c}
1\\
 -1\\
 \alpha\\
 \alpha
 \end{array}\right)$,
${\bf v}_2=\left(\begin{array}{c}
-1\\
1\\
\alpha\\
\alpha\end{array}\right)$, ${\bf v}_3=\left(\begin{array}{c}
 -1\\
 -1\\
 -\beta\\
 \beta\end{array}\right)$,
 ${\bf v}_4=\left(\begin{array}{c}
1\\
1\\
-\beta\\
\beta\end{array}\right)$
respectively. So the linear
transformation $x=Py$, where $P=\left(\begin{array}{cccc}
1& -1&-1&1\\
-1&1&-1&1\\
\alpha& \alpha& -\beta& -\beta\\
\alpha&\alpha&\beta&\beta\end{array}\right)$, turn the equation
(1.3) into
\begin{equation}
\dot {\bf y}(t)=B{\bf y}+{\bf f}(t),
\end{equation}
where ${\bf y}=\left(
\begin{array}{c}
y_1 \\
y_2\\
y_3\\
y_4\end{array}\right)$,
$B=
\begin{pmatrix}
\alpha & 0&0&0\\
0 &-\alpha&0&0\\
0&0&\beta&0\\
0&0&0&-\beta\end{pmatrix}$ and ${\bf f}(t)=P{\bf g}(t)=
\begin{pmatrix}
-g(t)-g(-t)\\
-g(t)-g(-t)\\
\beta(-g(t)+g(-t))\\
\beta(g(t)-g(-t))\end{pmatrix}. $

\medskip
\section{Case 1: $a < b < -a$}
By the standard formulation we can obtain following lemma.

{\lem Suppose that $g(t)\in C({\bf R})$, and bounded on ${\bf R}$,
and $\alpha>0,\beta>0$ .Then the general solution of the system
(1.4) on ${\bf R}$, is given by
\begin{equation}
{\bf y}=e^{{B}t}{\bf c}+ Y(t),
\end{equation}
where ${\bf c}=\begin{pmatrix}
c_1\\
c_2\\
c_3\\
c_4\end{pmatrix}$, and $Y(t)=\begin{pmatrix}
  \int_t^{\infty} e^{\alpha(t-s)}(g(s)+g(-s))ds\\
  -\int_{-\infty}^te^{-\alpha(t-s)}(g(s)+g(-s))ds\\
  -\beta\int_t^{\infty}e^{\beta(t-s)}(-g(s)+g(-s))ds\\
 \beta\int_{-\infty}^t e^{-\beta(t-s)}(g(s)-g(-s))ds
\end{pmatrix}$ and so the general solution of the system (1.2) or
(1.3) on ${\bf R}$, is given by
\begin{equation}{}
{\bf x}=P{\bf y}=P(e^{Bt}{\bf c}+ Y(t))
\end{equation}
or

\begin{equation}
x_1=c_1e^{\alpha t}-c_2e^{-\alpha t}-c_3e^{\beta t}+c_4e^{-\beta
t}+w_1(t)
\end{equation}
\begin{equation}
x_2=-c_1e^{\alpha t}+c_2e^{-\alpha t}-c_3e^{\beta t}+c_4e^{-\beta
t}+w_2(t)
\end{equation}
\begin{equation}
x_3=-c_1\alpha e^{\alpha t}+c_2\alpha e^{-\alpha t}-c_3\beta
e^{\beta t}+c_4\beta e^{-\beta t}+w_3(t)
\end{equation}
\begin{equation}
x_4=-c_1\alpha e^{\alpha t}+c_2\alpha e^{-\alpha t}+c_3\beta
e^{\beta t}+c_4\beta e^{-\beta t}+w_4(t)
\end{equation}
where
\begin{eqnarray*}
w_1(t)&=&\int_t^{\infty}
e^{\alpha(t-s)}(g(s)+g(-s))ds\\
&+&\int_{-\infty}^te^{-\alpha(t-s)}(g(s)+g(-s))ds\\
&+&\beta\int_t^{\infty}e^{\beta(t-s)}(-g(s)+g(-s))ds\\
&+&\beta\int_{-\infty}^te^{-\beta(t-s)}(g(s)-g(-s))ds.
\end{eqnarray*}
\begin{eqnarray*}
w_2(t)&=&\int_t^{\infty}
e^{\alpha(t-s)}(-g(s)-g(-s))ds\\
&+&\int_{-\infty}^te^{-\alpha(t-s)}(-g(s)-g(-s))ds\\
&+&\beta\int_t^{\infty}e^{\beta(t-s)}(-g(s)+g(-s))ds\\
&+&\beta\int_{-\infty}^te^{-\beta(t-s)}(g(s)-g(-s))ds.
\end{eqnarray*}

\begin{eqnarray*}
w_3(t)&=&\alpha\int_t^{\infty}
e^{\alpha(t-s)}(g(s)-g(-s))ds\\
&-&\alpha\int_{-\infty}^te^{-\alpha(t-s)}(-g(s)-g(-s))ds\\
&+&\beta^2\int_t^{\infty}e^{\beta(t-s)}(-g(s)+g(-s))ds\\
&-&\beta^2\int_{-\infty}^te^{-\beta(t-s)}(g(s)-g(-s))ds.
\end{eqnarray*}

\begin{eqnarray*}
w_4(t)&=&\alpha\int_t^{\infty}
e^{\alpha(t-s)}(g(s)+g(-s))ds\\
&-&\int_{-\infty}^te^{-\alpha(t-s)}(g(s)+g(-s))ds\\
&-&\beta^2\int_t^{\infty}e^{\beta(t-s)}(-g(s)+g(-s))ds\\
&+&\beta^2\int_{-\infty}^te^{-\beta(t-s)}(g(s)-g(-s))ds.
\end{eqnarray*}

}

Because there are four arbitrary constants in (3.3), we can not
conclude that (3.2) is the general solution of (1.1). (3.3) may
not even be a solution of (1.1) for some constants $c_1,c_2,c_3$
and $c_4$ indeed.

{\lem  Let $g(t)\in C({\bf R})$ is bounded on ${\bf R}$, and
$\alpha>0,\beta>0$, then every solution of Eq.(1.1) is of the form
(3.3), if and only if $c_1=c_2, c_3=-c_4$, that is the general
solution of Eq.(1.1) is of the form
\begin{equation}\label{}
x(t)=k_1(e^{\alpha t}-e^{-\alpha t})+k_2(e^{\beta t}+e^{-\beta
t})+w_1(t)
\end{equation} where $k_1,k_2$ are arbitrary constants. }

{\it Proof}\,\, From the requirements of $x_1=x(t)$,
$x_2=x_1(-t)$, $x_3=\dot x_1(t)$, $x_4=x_3(-t)$, we can derive
$c_1=c_2, c_3=-c_4$. Let $k_1=c_1=c_2, k_2=-c_3=c_4$.

{\thm Under the condition of Lemma (3.1), Eq.(1.1) has a unique
bounded solution given by
$$
\begin{array}{c} x(t)=w_1(t)=\int_t^{\infty}
e^{\alpha(t-s)}(g(s)+g(-s))ds+\int_{-\infty}^te^{-\alpha(t-s)}(g(s)+g(-s))ds\\
+\beta\int_t^{\infty}e^{\beta(t-s)}(-g(s)+g(-s))ds+\beta\int_{-\infty}^te^{-\beta(t-s)}(g(s)-g(-s))ds.
\end{array}
$$ and moreover, $\sup_{t\in
R}|x(t)|\leq (\frac{2}{\alpha}+1) \sup_{t\in R}|g(t)|$}

{\it Proof}\,\, (3.11) is a bounded solution if and only if
$k_1=k_2=0 $. We obtain the inequality $\sup_{t\in R}|x(t)|\leq
(\frac{2}{\alpha}+1) \sup_{t\in R}|g(t)|$ by evaluating
$x(t)=w_1(t)$.

 {\thm Let $g(t)\in AP({\bf R})$, then Eq.(1.1) has a unique
almost periodic solution $x(t)$, and $mod(x)=mod(g)$.
Furthermore,if $g(t)$ is periodic, then Eq.(1.1) has a unique
harmonic solution.}

{\it Proof}\,\, We will show $x(t)=w_1(t)$ is almost periodic
solution of (1.1). For $\tau \in T(g(t),\epsilon)$,
\begin{eqnarray*}
&&|x(t+\tau)-x(t)|\\
&=&\left|\int_{t+\tau}^{\infty}e^{\alpha(t+\tau-s)}(g(s)+g(-s))ds-\int_{-\infty}^{t+\tau}e^{-\alpha(t+\tau-s)}(-g(s)-g(-s))ds\right.\\
&&+\beta\int_{t+\tau}^{\infty}e^{\beta(t+\tau-s)}(-g(s)+g(-s))ds+\beta\int_{-\infty}^{t+\tau}e^{-\beta(t+\tau-s)}(g(s)-g(-s))ds\\
&-&\int_t^{\infty}e^{\alpha(t-s)}(g(s)+g(-s))ds+\int_{-\infty}^te^{-\alpha(t-s)}(-g(s)-g(-s))ds\\
&&-\left.\beta\int_t^{\infty}e^{\beta(t-s)}(-g(s)+g(-s))ds-\beta\int_{-\infty}^te^{-\beta(t-s)}(g(s)-g(-s))ds\right|\\
&\leq&\int_t^{\infty}e^{\alpha(t-s)}(|g(s+\tau)-g(s)|+|g(-s-\tau)-g(-s)|)ds\\
&&+\int_{-\infty}^te^{-\alpha(t-s)}(|g(s+\tau)-g(s)|+|g(-s-\tau)-g(-s)|)ds\\
&&+\beta\int_t^{\infty}e^{\beta(t-s)}(|g(s+\tau)-g(s)|+|g(-s-\tau)-g(-s)|)ds\\
&&+\beta\int_{-\infty}^te^{-\beta(t-s)}(|g(s+\tau)-g(s)|+|g(-s-\tau)-g(-s)|)ds\\
&\leq&4(\frac{1}{\alpha}+1)\epsilon.
\end{eqnarray*}
So $x(t)\in AP({\bf R})$, and $mod(x)\subset mod(g)$. From
$g(t)=\ddot{x}(t)+ax(t)+bx(-t)$, and the lemma 2.4 we conclude
$mod(g)\subset mod(x)$, and so $mod(x)= mod(g)$. If $g(t)$ is
$\omega-$periodic, then $|x(t+\omega)-x(t)|\leq
4(\frac{1}{\alpha}+1) ||g(t+\omega)-g(t)||=0$. So $x(t)$ is
$\omega-$periodic solution.

Uniqueness. If there is another almost periodic solution $x_1(t)$
for Eq.(3), then the difference $x(t)-x_1(t)$ should be a solution
of the homogeneous equation
\begin{equation}\label{}
\ddot {x}(t)+ax(t)+bx(-t)=0,b\neq 0, t\in {\bf R}.
\end{equation}
According to the Lemma 3.12, we can derive
\begin{equation}\label{}
 x(t)-x_1(t)=k_1(e^{\alpha t}-e^{-\alpha t})+k_2(e^{\beta t}-e^{-\beta
 t})
\end{equation}
for some constant $k_1,k_2$. If $|k_1|+|k_2|\neq 0$ , then $
x(t)-x_1(t)$ will be unbounded. This is a contradiction to the
boundedness of almost periodic function. So $ x(t)-x_1(t)\equiv
0$, i.e $ x(t)\equiv x_1(t)$.

{\thm Suppose $f(t,x,y)$ is almost periodic on t uniformly with
respect to $x$ and $y$ on any compact set ${\bf W} \subset {\bf C
}^2$, and satisfies Lipschitz condition
$$
|f(t,x_1,y_1)-f(t,x_2,y_2)|\leq L(|x_1-x_2|+|y_1-y_2|)
$$
for any $(x_1,y_2),(x_2,y_2)\in {\bf W}$, where
$L<\frac{\alpha}{4(1+\alpha)}$. Then Eq.(1.2) has a unique almost
periodic solution $x(t)$ and $mod(x)=mod(f).$  In addition, if $f$
 is periodic in $t$, then Eq.(1.2) has a unique harmonic solution
$x(t)$. }

{\it Proof}\,\,\,\,We know the subset
$${\bf B}=\{\phi(t):\phi \in AP({\bf R}),mod(\phi)\subset mod(f) \}$$
 of $AP({\bf R})$ is a Banach space with the supremum norm $||\phi||=\sup_{t\in R}|\phi(t)|$.
For any $\phi \in {\bf B}$, we know $f(t,\phi(t),\phi(-t))\in {\bf
B}$ [11]. According to the theorem 2.1 we see the equation
\begin{equation}\label{}
  \ddot {x}(t)+ax(t)+bx(-t)=f(t,\phi(t),\phi(-t)),b\neq 0,t\in R
\end{equation}
possess a unique almost periodic solution, denote it by
$(T\phi)(t)$. Then we define a mapping $T: {\bf B}\rightarrow {\bf
B}$. Now we show $T$ is contracted.

For $\phi(t),\psi(t)\in {\bf B}$, the equation
\begin{equation}\label{}
\ddot{x}(t)+ax(t)+bx(-t)=f(t,\phi(t),\phi(-t))-f(t,\psi(t),\psi(-t)),
b\neq 0,t\in R
\end{equation}
has a unique almost periodic solution $(T\phi-T\psi)(t)$, and
\begin{eqnarray*}
&&|(T\phi-T\psi)(t)|\\
&=&|\int_t^{\infty}e^{\alpha(t-s)}[f(s,\phi(s),\phi(-s))-f(s,\psi(s),\psi(-s))\\
&&+f(-s,\phi(-s),\phi(s))-f(-s,\psi(-s),\psi(s))]ds\\
&+&\int_{-\infty}^te^{-\alpha(t-s)}[f(s,\phi(s),\phi(-s))-f(s,\psi(s),\psi(-s))\\
&&+f(-s,\phi(-s),\phi(s))-f(-s,\psi(-s),\psi(s))]ds\\
&+&\beta\int_t^{\infty}e^{\beta(t-s)}[-f(s,\phi(s),\phi(-s))-f(s,\psi(s),\psi(-s))\\
&&+f(-s,\phi(-s),\phi(s))-f(-s,\psi(-s),\psi(s))]ds\\
&+&\beta\int_{-\infty}^te^{-\beta(t-s)}[f(s,\phi(s),\phi(-s))-f(s,\psi(s),\psi(-s))\\
&&-f(-s,\phi(-s),\phi(s))-f(-s,\psi(-s),\psi(s))]ds|\\
&\leq& 4L(\frac{1}{\alpha}+1)\|\phi-\psi\|.
\end{eqnarray*}
So
$$
\parallel T\phi-T\psi \parallel \leq
4L(\frac{1}{\alpha}+1)\|\phi-\psi\|.
$$

Since $L<\frac{\alpha}{4(1+\alpha)}$,  $T$ is a contraction
mapping, and so $T$ has a unique fixed point in $B$. That is to
say the equation (1.1) has a unique almost periodic solution
$x(t)$ and $mod(x)\subset mod(f)$.

 If $f(t,x,y)$ is continuous $\omega$-periodic in $t$, then for any
 $\omega$-periodic function $\phi (t)$, $f(t,\phi (t),\phi(-t))$
 is continuous $\omega$-periodic function too. We denote by $C_{\omega}$ the set of all continuous
 $\omega$-periodic functions. Then we know $C_{\omega}$ is a Banach space with
 supremum norm $||\phi||=\sup_{t\in R}|\phi(t)|$. From theorem 2.1,
 we conclude, for any $\phi(t)\in C_w$, Eq.(2.4) has a unique $\omega$-periodic
 solution $T\phi\in C_w$. We can easily prove as above that $T$ is
 contracted. So $T$ has a unique fixed point $x(t)\in C_w$, i.e.
 there is a unique harmonic solution for Eq.(1.1). So the theorem
 2.2 is completed.
 \medskip

 \section{Case 2: $-a < b < a $}

 In this case, both $\alpha$ and $\beta$ are pure imaginary
 numbers. Set $\alpha=i \mu $ and $\beta =i \nu $. By the standard formulation we can obtain

 {\lem Suppose that $g(t)\in C({\bf R})$, and
$-a < b < a $. Then the general solution of the system (2.3) on
$R$, is given by
\begin{equation}
{\bf y}=e^{{B}t}\left({\bf c}+ \int_0^te^{{-B}s}{\bf f}(s)ds
\right),
\end{equation}
and so the general solution of (2.2) is given by
\begin{equation}{}
{\bf x}=P{\bf y}=Pe^{Bt}\left({\bf c}+ \int_0^te^{{-\bf B}s}{\bf
f}(s)ds \right)
\end{equation}
where ${B}=diag\{i\mu ,-i\mu ,i\nu , -i\nu  \} $ }.

Similar to section 3 we can derive that the general solution of
Eq.(1.2) take the form

\begin{equation}
x(t)=k_1(e^{i\mu t}-e^{-i\mu t})+k_2(e^{i\nu t}+e^{-i\nu t})
    -e^{i\mu t}\int_0^te^{-i\mu s}(g(s)+g(-s))ds
\end{equation}
$$+e^{-i\mu t}\int_0^te^{i\mu
  s}(g(s)+g(-s))ds  -\beta e^{i\nu t}\int_0^te^{-i\nu s}(-g(s)+g(-s))ds
$$
$$+\beta e^{-i\nu t}\int_0^te^{i\nu s}(g(s)-g(-s))ds
$$

To this end,we introduce a result of Favard.
 {\lem$^{[8]}$Assume the fourier series of almost function $f(t)$ is $\Sigma A_ne^{i\Lambda_nt}$.If
there is $\alpha > 0$ such that $\Lambda_n>\alpha$, for any $n\in
Z$, then the indefinite integral of $f(t)$ is almost periodic }

Then we have following theorem.

{\thm If $\forall \lambda_k\in Freq(g)$ s.t $|\mu \pm
\lambda_k|\geq \rho>0$, $|\nu \pm \lambda_k|\geq \sigma>0$ and $-a
< b < a $, where $\rho$ and $\sigma$ are some constants, then (1)
every solution of (1.1) is almost periodic; (2) there exists a
unique almost periodic solution $x(t)$ satisfies initial condition
$x(0)=x_0$, $\dot x(0)=\dot x_0$ }

{\it Proof}\,\,\, (1)  Lemma 2.4 and Lemma 4.2 yield those
indefinite integrals in $(4.3)$ are almost periodic.

(2) A solution $x(t)$ s.t $x(0)=x_0$,$\dot x(0)=\dot x_0$ if and
only if $k_1=\dot x_0/2i\mu$ and $k_2=x_0/2$ in $(4.3)$.

{\rem \,\,The existence of the term $bx(-t)$ with reflection of
the argument in (1.1) may influence the boundeness of the
solutions of the equation without reflection of the argument
drastically.

 Example 1 \,\ The general solution of second order equation
\begin{equation}
\ddot x+4x=\cos2t
\end{equation}
is
$$
x(t)=c_1\cos2t+c_2\sin2t+\frac{t}{4}\sin2t.
$$
So the equation (4.3) has no bounded solution. But according to
Theorem 4.1, every solution of the equation
$$
\ddot x(t)+4x(t)+2x(-t)=\cos2t
$$
 is bounded, and almost periodic.

 Example 2 \,\,The general solution of the equation
 \begin{equation}
 \dot x(t)+2x(t)=\cos t
 \end{equation}
 is $x(t)=c_1\cos\sqrt{2}t+c_2\sin\sqrt{2}t+\cos t$, so every solution of (4.5) is bounded. But from formulae (4.3),
 we can derive every solution of
 $$
  \dot x(t)+2x(t)+x(-t)=\cos t
 $$ is unbounded.  We note that $\{-1, 1\} \subset Freq(g)=Freq(\cos t)$, and
 $\mu=\sqrt{a-b}=1$. This shows the condition $\rho>0, \sigma>0$ of Theorem 4.2 is
 sharp.}

{\rem \,\, One can investigate the solution of (1.1) and (1.2) for
the other case of $a$ and $b$. It seems that the study for the
case $a=b$ or $a=-b$ is complicated. To study the solutions for
nonlinear equation
$$
x''+x^3+x(-t)=p(t)
$$
may be vary interesting problem.}

\bibliographystyle{amsplain}

\bibliography{}
\end{document}